\documentclass[a4paper,reqno,11pt]{amsart}
\usepackage{amssymb,amsmath,amsopn}
\usepackage{graphics} %This is used to rotate labels in figure 1
\usepackage{setspace} %Used for submitted draft

\newtheorem{theorem}{Theorem}[section]

\newtheorem{lemma}[theorem]{Lemma}
\newtheorem{corollary}[theorem]{Corollary}
\newtheorem{proposition}[theorem]{Proposition}
\newtheorem{question}[theorem]{Question}
\newtheorem{problem}[theorem]{Problem}

\DeclareMathOperator{\sgn}{sgn}
\DeclareMathOperator{\supp}{supp}

\newcommand{\N}{\mathbf{N}}
\newcommand{\Z}{\mathbf{Z}}
\newcommand{\Q}{\mathbf{Q}}
\newcommand{\R}{\mathbf{R}}

\newcommand{\Ind}{\big\uparrow}
\newcommand{\Res}{\big\downarrow}
\newcommand{\ind}{\!\!\uparrow}
\newcommand{\res}{\!\!\downarrow}

\newcommand{\mfrac}[2]{ { \textstyle{\frac{#1}{#2}} } }

\newcounter{thmlistcnt}
\newenvironment{thmlist}%
	{\setcounter{thmlistcnt}{0}%
	\begin{list}{\emph{(\roman{thmlistcnt})}}{%
		\usecounter{thmlistcnt}%
		\setlength{\topsep}{0pt}%
		\setlength{\leftmargin}{0pt}%
		\setlength{\itemsep}{0pt}%
		\setlength{\itemindent}{25pt}}%
	}%
	{\end{list}}%

\newcounter{indentedthmlistcnt}
\newenvironment{indentedthmlist}%
	{\setcounter{indentedthmlistcnt}{0}%
	\begin{list}{\emph{(\roman{indentedthmlistcnt})}}{%
		\usecounter{indentedthmlistcnt}%
		\setlength{\labelwidth}{32pt}%
		\setlength{\topsep}{-3pt}%
		\setlength{\leftmargin}{32pt}%
		\setlength{\itemsep}{0pt}%
		\setlength{\itemindent}{0pt}}%
	}%
	{\end{list}}%

\newcounter{rmklistcnt}
\newenvironment{rmklist}%
	{\setcounter{rmklistcnt}{0}%
	\begin{list}{(\arabic{rmklistcnt})}{%
		\usecounter{rmklistcnt}%
		\setlength{\topsep}{0pt}%
		\setlength{\leftmargin}{0pt}%
		\setlength{\itemsep}{3pt}%
		\setlength{\itemindent}{0pt}
		\setlength{\leftmargin}{19pt}} %18pt}}%
	}%
	{\end{list}}%
	
\newcounter{indentedrmklistcnt}
\newenvironment{indentedrmklist}%
	{\setcounter{indentedrmklistcnt}{0}%
	\begin{list}{(\arabic{indentedrmklistcnt})}{%
		\usecounter{indentedrmklistcnt}%
		\setlength{\topsep}{6pt}%
		\setlength{\leftmargin}{0pt}%
		\setlength{\itemsep}{6pt}%
		\setlength{\itemindent}{31pt}
		\setlength{\leftmargin}{0pt}} %18pt}}%
	}%
	{\end{list}}%

\begin{document}

\begin{abstract}
By exploiting relationships between the values taken by
ordinary characters of  symmetric groups we prove
two theorems in the modular representation theory of the 
symmetric group. 

1. The decomposition matrices of symmetric groups in odd 
characteristic have distinct rows. In characteristic $2$ the rows 
of a decomposition matrix labelled
by the different partitions~$\lambda$ and~$\mu$ are equal if and only 
if~$\lambda$ and~$\mu$ are conjugate. An analogous result is proved
for Hecke algebras.

2. A Specht module for the symmetric group $S_n$, defined over an 
algebraically closed field of odd characteristic, is decomposable on 
restriction to the alternating group $A_n$ if and only if it is simple, 
and the labelling partition is self-conjugate. This result is generalised
to an arbitrary field of odd  characteristic.
\end{abstract}

%\onehalfspace

\title[Decomposition matrices of symmetric groups]
{Character values and decomposition matrices of symmetric groups}

\author{Mark Wildon}
\date{October 12, 2006 } 
%Used in LMS: Author's email address: \texttt{wildon@maths.ox.ac.uk}}
\email{m.j.wildon@swansea.ac.uk}
%Used in LMS: \classno{20C30 (primary), 20C20 (secondary)}

\maketitle
\thispagestyle{empty}

\section{Introduction}
In this paper we solve two problems in the modular 
representation theory of the symmetric group. The 
first asks for a necessary and sufficient condition
for two rows of a decomposition matrix of a symmetric
group to be equal. The second asks for a characterisation of
the Specht modules which decompose on restriction
from the symmetric group to the alternating group.
Although these problems may seem quite different
from one another, both can be solved 
by  similar arguments using the ordinary characters of
the symmetric group. In fact, both problems can be reduced
to questions typified by the following:

\begin{question}
Suppose that two ordinary irreducible characters of the symmetric group agree
on all elements of order not divisible by $3$ (that is, $3'$-elements)
--- must they be the same?
\end{question}

We give a general strategy for answering questions such as this in~\S 2 below. 
Our idea is to use the central characters of symmetric groups 
to find algebraic relationships between the values taken by a fixed ordinary 
irreducible character on different conjugacy classes. The main results we 
prove may be found below in Corollaries~\ref{Cor:Prop},~\ref{Cor2:Prop}, 
and~\ref{Cor:Prop2}.

To give a representative example, Corollary~\ref{Cor2:Prop} implies that, 
given the values taken by an ordinary irreducible character of a symmetric group
on~$3'$-elements, one can determine all its remaining values. Thus the
question posed above has an affirmative answer. As this example may suggest, 
our results on character values are of some independent interest. In~\S 2.5 
we give some questions they inspire.

We now outline the problems that will be solved using the results of \S 2.

\subsection{Decomposition matrices}
A \emph{partition} of a number $n \in \N$ is a sequence 
$(\lambda_1, \ldots, \lambda_k)$ of positive integers such that 
$\lambda_1 \ge \lambda_2 \ge \ldots \ge \lambda_k \ge 1$
and $\lambda_1 + \ldots + \lambda_k = n$. To indicate that~$\lambda$ is 
a partition of~$n$ we write~$\lambda \vdash n$.

Let $F$ be a field and let~$S^\lambda$ be the Specht module for $FS_n$
labelled by the partition~$\lambda$ of $n$. For the definition and some
examples of these modules see Chapters~4 and~5 of~\cite{James}. 
We recall here that if~$F$ has characteristic zero then every Specht module 
is simple, and every simple~$FS_n$-module is isomorphic to a Specht module.
If~$F$ has prime characteristic~$p$ then this is no longer the case. 
However, if~$\lambda$ is $p$-regular --- that is,~$\lambda$ has no more 
than $p-1$ parts of any given size --- then~$S^\lambda$ has a simple top,
denoted~$D^\lambda$. The modules~$D^\lambda$ are pairwise
non-isomorphic and give all the simple representations of~$FS_n$.
We record the composition factors of Specht modules in 
characteristic~$p$ in the \emph{decomposition matrix}~$D_p(n)$, 
defined by letting~$D_p(n)_{\lambda \nu}$ be the number of composition 
factors of~$S^\lambda$ that are isomorphic to~$D^\nu$.

A fundamental problem in modular representation
theory is to determine the decomposition matrices of symmetric groups. 
In~\S 3 we prove the following theorem.

\begin{theorem}\label{Thm:DecDistinct}
Let $p$ be prime and let $n \in \N$.

\begin{indentedthmlist}
\item If $p > 2$ then the rows of $D_p(n)$ are mutually distinct.

\item If $p=2$ then the rows labelled by $\lambda$ and $\mu$ are the same
if and only if $\lambda = \mu$ or $\lambda = \mu'$, the conjugate
partition to $\mu$.
\end{indentedthmlist}
\end{theorem}

\noindent Thus in odd characteristic, a Specht module is determined
by its set of composition factors. In characteristic $2$, there
are at most two Specht modules with any given set of composition
factors. (For the definition of the conjugate of 
a partition  see~\cite[Definition~3.5]{James}.)

\medskip
\noindent\emph{Remarks on Theorem \emph{\ref{Thm:DecDistinct}}.}

\begin{rmklist}
\item
It is well known (see \cite[Corollary 12.3]{James}) 
that when the partitions labelling the
rows and columns of a decomposition matrix
are ordered lexicographically, but with $p$-regular
partitions placed before non-$p$-regular partitions, 
the matrix takes a `wedge' shape, illustrated below by $D_3(5)$.

\begin{figure}[h]
\newcommand{\rb}[1]{\rotatebox{90}{#1}}
\newcommand{\cd}{\cdot}
\[ \begin{matrix} 
	& \rb{(5)} & \rb{(4,1)} & \rb{(3,2)} & \rb{(3,1,1)} & \rb{(2,2,1)}  \\
	\hfill (5)    & 1 & \cd & \cd & \cd & \cd \\
	\hfill(4,1)   & \cd & 1 & \cd & \cd & \cd \\
	\hfill(3,2)   & \cd & 1 & 1 & \cd & \cd \\
	\hfill(3,1,1) & \cd & \cd & \cd & 1 & \cd \\
	\hfill(2,2,1) & 1 & \cd & \cd & \cd & 1 \\
	\hfill(2,1,1,1)&\cd & \cd & \cd & \cd & 1 \\ 
	\hfill(1^5)   & \cd & \cd & 1 & \cd & \cd 
\end{matrix} \]
\caption{The decomposition matrix of $S_5$ in characteristic $3$.}
\end{figure}

It is therefore easy to distinguish between the rows labelled by
$p$-regular partitions. The
force of Theorem~\ref{Thm:DecDistinct} comes from the
fact that, when~$n$ is large compared to~$p$,
most partitions are not $p$-regular. More precisely,
if for~\hbox{$\ell \ge 2$} we let $g_\ell(n)$ be the proportion
of~$\ell$-regular partitions of~$n$
(here~$\ell$ is not necessarily prime) then
\[ g_\ell(n) \sim A n^{1/4} \thinspace \mathrm{e}^{-c 
\left( 1- \sqrt{\frac{\ell-1}{\ell}} \right) \sqrt{n}} 
\quad\text{as $n\rightarrow \infty$} \]
where $c = 2\sqrt{\pi^2/6}$ and  $A \in \R$ depends only on $\ell$. 
The proportion of $\ell$-regular partitions therefore tends rapidly to zero. 
This formula was  proved by Hagis using the circle-method
(see \cite[Corollary~4.2]{Hagis}). It is  interesting  to see how close one
can get to it by less sophisticated methods.
When $\ell = 2$ I have given an elementary proof 
(see \cite[\S 5]{WildonAbacus}), but when $\ell > 2$, the
strongest result I have been able to obtain by elementary methods is
\[ \log g_\ell(n) 
\sim -2\sqrt{\frac{\pi^2}{6}} 
\left( 1- \sqrt{\frac{\ell-1}{\ell}} \right) \sqrt{n} < -\frac{1}{\ell}\sqrt{\frac{\pi^2}{6}} \sqrt{n}. 
\]

\item The analogue of Theorem~\ref{Thm:DecDistinct} 
for the Hecke algebras of symmetric groups may also be
proved using the results of \S 2 --- see Theorem~\ref{Thm:DecDistinctHecke}.
I hope to report later on the situation for alternating groups; 
for a partial result see Theorem~\ref{Thm:Alt}. 
As Schur algebras have lower-unitriangular decomposition
matrices (see \cite[Theorem 3.5a]{GreenGLn}),
the rows of their decomposition matrices are always distinct. 

\item When defined over a field of characteristic $2$, Specht modules
labelled by different partitions may be isomorphic.
Theorem \ref{Thm:DecDistinct}(ii) can be used to show
that~$S^\lambda$ is isomorphic to $S^\mu$ if and only either
$\lambda = \mu$, or $\lambda = \mu^\prime$ and~$S^\lambda$ is self-dual. 
Unfortunately it does not seem easy
to classify the self-dual Specht modules in characteristic~$2$. 
For example, in characteristic $2$, $S^{(5,2)}$ is simple, and hence 
self-dual, but~$S^{(5,1,1)}$ is also self-dual, and even decomposable 
(see~\cite[\S 23.10]{James}).
\end{rmklist}

\subsection{The restriction of Specht modules to the alternating group}

Our main result is the following theorem, which we prove in \S 4.

\begin{theorem}\label{Thm:SpechtRes}
Let $F$ be an algebraically closed field of odd characteristic. 
Let~$\lambda$ be a partition of $n\in \N$. The Specht module
 $S^\lambda$ defined over $F$ is decomposable on restriction to $A_n$ if and only 
if $S^\lambda$ is simple and $\lambda = \lambda^\prime$.
\end{theorem}

\noindent\emph{Remarks on Theorem \emph{\ref{Thm:SpechtRes}}.}

\begin{rmklist}
\item
A related result is proved in \cite{FordIrreps}, where Ford uses the
Ford--Kleshchev  proof \cite{FordKleshchev} of the Mullineux conjecture 
to give a straightforward way of determining the simple
$FS_n$-modules $D^\lambda$ such that $D^\lambda \cong D^\lambda \otimes \sgn$,
and hence (using some basic Clifford theory which we repeat in \S 4) 
those irreducible representations
of $FS_n$ which split on restriction to~$A_n$. 

\item
In \cite[Conjecture~5.47]{MathasHecke} James and Mathas conjectured a 
necessary and sufficient condition for a Specht module to be simple.
Their conjecture was subsequently proved by Fayers (see
\cite{FayersIf, FayersOnlyIf}). His work makes
it a simple matter to work with the criterion given in our theorem. 
In \S 4.1 we state the James--Mathas condition and use
it to generalise Theorem~\ref{Thm:SpechtRes}
to fields that are not algebraically closed. 
\end{rmklist}

\section{Results on character values}
Our approach uses the central characters of symmetric groups and some 
elementary properties of the sums of elements in symmetric group 
conjugacy classes.
We work all the time over the field of rational numbers.

First we must
introduce some notation. Fix $n \in \N$.
For $i \in \N$ we let~$s_i \in \Q S_n$ be the sum of all~$i$-cycles 
in~$S_n$,
so in particular $s_1 = 1_{S_n}$, the identity element in~$\Q S_n$.
More generally, if~$\mu$ is a partition of $n$, 
we let~$s_\mu$ be the sum of all elements in $S_n$
of cycle type~$\mu$. When writing the elements~$s_\mu$ we shall simplify
the notation by always
ignoring parts of size $1$; thus~$s_{(i)}$ is the same as~$s_i$, and if, for example $n=9$,
then $s_{(3,2)} = s_{(3,2,1^4)}$. (This leads to no ambiguity,
as the degree~$n$ is fixed throughout this section.) If $\mu$
has $m$ parts of size $1$ then we say that $\mu$ has
\emph{support}~$n-m$, and write~$\supp \mu = n-m$.
Let~$K_\mu$ be the number
of elements of~$S_n$ of cycle type~$\mu$. 

If $\lambda$ is a partition of $n$, we write~$\chi^\lambda$ for the
the irreducible ordinary character of~$S_n$ afforded by
the Specht module~$S^\lambda$ (now defined over a field of characteristic 
zero). Let~$\chi^\lambda(\mu)$ be the value of~$\chi^\lambda$ on elements of cycle type~$\mu$. 
The \emph{central character} corresponding to~$\lambda$ is the algebra
homomorphism $\omega_\lambda : Z(\Q S_n) \rightarrow \Q $ defined 
by mapping an element~$z \in Z(\Q S_{n})$ to the scalar by which it 
acts on the Specht module~$S^\lambda$. As
\begin{equation}\label{eq:cent} 
\omega_\lambda (s_\mu)  = \frac{\chi^\lambda(\mu)K_\mu}{\chi^\lambda(1)}
\end{equation}
the values of $\chi^\lambda$ are determined by  $\omega_\lambda$.
Since the $s_\mu$ for $\mu \vdash n$ form a linear
basis for $Z(\Q S_n)$, the converse also holds.

\subsection{Generating sets for $Z(\Q S_n)$} 
It appears to have first been proved by Kramer \cite{Kramer} that
the cycle sums $s_1, \ldots, s_n$ generate $Z(\Q S_n)$ as a 
$\Q$-algebra. It therefore follows from \eqref{eq:cent} 
that given any partition $\mu$ of $n$, there
is a polynomial $P_\mu(X_1,\ldots,X_n) \in \Q [X_1,\ldots,X_n]$ such that
\[ \chi^\lambda(\mu) = P_\mu\bigl(\chi^\lambda(1_{S_n}),\chi^\lambda((1 2)),
\ldots, \chi^\lambda((1\,2\ldots n))\bigr) \]
where $(1\,2\ldots i)$ is the standard $i$-cycle in~$S_n$.
It is an interesting feature of these polynomials that
they are entirely independent of the partition~$\lambda$. 

As an example, we determine $P_{(2,2,1^{n-4})}$. 
A short calculation shows that in~$\Q S_n$,
$s_{(2)}^2 = 2s_{(2,2)} + 3s_{(3)} + n(n-1)/2$, and hence
\[ P_{(2,2,1^{(n-4)})} = -\frac{n(n-1)}{4} X_1 + \frac{1}{2} X_2^2 - 
\frac{3}{2} X_3.\]
Another consequence of this relation is that 
$\chi^\lambda((123))$ is determined by the values
of $\chi^\lambda$ at the $3'$-elements $1_{S_n}$,~$(12)$, and $(12)(34)$. 
A straightforward generalisation of this, given in part (ii) 
of the proposition below, can be used to answer Question~1.1. 

In connection with these polynomials, it is worth mentioning that the values
taken by ordinary characters on cycles may easily be calculated  using
the Murnagham--Nakayama rule (see \cite[Ch.~21]{James}).
Another way to find these values, of more theoretical
interest, is to use the combinatorial interpretation of 
Frumkin, James and Roichman \cite{FrumkinCycles}.
For cycles of small length there are also some
interesting explicit formulae (see for instance~\cite{Ingram}).

\begin{proposition}\label{Prop:ThmAProp}
Let $\mu$ be a partition of $n$ and let $\ell \in \N$.

\begin{thmlist}
\item 
The $\Q $-algebra generated by
\[ \left\{ s_i : i \le \supp \mu \right\} \] 
contains $s_\mu$. Moreover, if $\mu$ labels a conjugacy class of 
odd permutations, then~$s_\mu$ may be expressed as a polynomial 
in the $s_i$ in such a way that in every monomial term in the expression,
at least one class sum $s_{2j}$ appears. 

\item
If $\ell > 2$ then the conjugacy class sums
\[ X_\ell (n) =
\bigl\{ \hbox{$s_i : 1 \le i \le n$, $\ell \not\hspace{3pt}\mid\; i$} \bigr\} \cup
	 \bigl\{ \hbox{$s_{(\ell j-1,2)} : 1 < \ell j < n$} \bigr\}\]
generate $Z(\Q S_n)$ as a $\Q$-algebra.

\item 
If $\ell$ is odd and $\ell > 1$ then the conjugacy class sums
\[ Y_\ell(n) = 
\bigl\{ \hbox{$s_i : 1 \le i \le n$, if $i$ is
even then $\ell \not\hspace{3pt}\mid\; i$} \bigr\} \cup
	 \bigl\{ \hbox{$s_{(2\ell j-1,2)} : 1 < 2\ell j < n$} \big\}\]
generate $Z(\Q S_n)$ as a $\Q$-algebra. 
Moreover, if $\mu$ labels a conjugacy class of odd permutations then
$s_\mu$ may be expressed as a polynomial in the elements of~$Y_\ell(n)$ in
such a way that in every monomial term in the expression,
at least one class sum of a conjugacy class of
odd permutations appears.
\end{thmlist}
\end{proposition}

\begin{proof} \vbox{Suppose that $\mu = (n^{a_n}, \ldots, 2^{a_2}, 1^{a_1})$. 

(i) We work by induction on $\supp \mu$. We have
\[ s_1^{a_1} s_2^{a_2} \ldots s_n^{a_n} = \alpha s_\mu + y \]
where $\alpha$ is a strictly positive integer
and~$y$ is a weighted sum of conjugacy class sums $s_\lambda$ for 
partitions~$\lambda$ of
support at most $\supp \mu -1$. Moreover, if~$\mu$ labels 
a conjugacy class of odd permutations then so does every $\lambda$ which
appears in~$y$. So, by induction,~$y$ may
be written as a polynomial in the~$s_i$ of the required form.}

(ii) Given (i), it is sufficient to prove that the conjugacy class
sums $s_{\ell j}$ are in the $\Q$-algebra generated
by $X_\ell(n)$.  For this, it is sufficient to prove by induction on~$j$  
that if  $1 < \ell j \le n$ then the conjugacy class sum~$s_{\ell j}$ 
is in the~$\Q$-algebra generated by 
\[
\bigl\{ \hbox{$s_i : 1 \le i < \ell j$}
\bigr\} \cup \bigl\{ s_{(\ell j - 1,2)} \bigr\}.\]
(The last term above should be disregarded if $\ell j = n$.)
If~$\ell j < n$ then
\begin{equation}\label{eq:1} s_2 s_{\ell j-1}= \alpha s_{\ell j} +
\beta s_{(\ell j-1,2)} + \!\sum_{1 \le i \le \ell j/2 - 1} 
\gamma_i s_{(\ell j-1-i,i)} 
\end{equation}
for some coefficients $\alpha,\beta,\gamma_i$, about which all we 
need to know is the obvious fact that~\hbox{$\alpha > 0$}. 
Each element of $S_n$ appearing in the conjugacy class sums 
$s_{(\ell j-1-i,i)}$ fixes one more point than a $\ell j$-cycle, 
so by~(i), each $s_{(\ell j-1-i,i)}$ is in the $\Q$-algebra generated by
$\left\{ s_i : i < \ell j \right\}$. Hence $s_{\ell j}$
is in the~$\Q$-algebra generated by 
$\left\{ s_i : 1 \le i < \ell j \right\} \cup \{s_{(\ell j - 1,2)}\}$,
as required.

If $\ell j = n$ then
\begin{equation}\label{eq:2} s_2s_{n-1} = \alpha^\prime \hspace{-1pt} s_{n} + 
\!\sum_{1 \le i \le n/2-1} \gamma_i^\prime s_{(n-1-i,i)}
\end{equation}
for some further coefficients $\alpha', \gamma_i'$. Again it is obvious
that~$\alpha' > 0$, so the result follows 
in the same way as before.

(iii) We use the same strategy as in (ii).
By (i) it suffices to prove that if
$1 < 2\ell j \le n$ then~$s_{2\ell j}$ may
be expressed as a polynomial in the elements
\[ \bigl\{ \hbox{$s_i : 1 \le i < 2\ell j$}
\bigr\} \cup \bigl\{ s_{(2\ell j - 1,2)} \bigr\}\]
in such a way that in every monomial in the expression, at least one 
class sum of odd permutations appears. As in (ii), 
this follows by induction on $j$, using 
equations~\eqref{eq:1}
and~\eqref{eq:2} above.
\end{proof}

In order to move from conjugacy class sums to individual
elements of~$S_n$ we introduce the following strategically chosen family of 
sets. For $\ell \in \N$, let
\[ \begin{split} Z_\ell(n) = \left\{ 
	(1\:2\ldots 2j) : \hbox{$1 \le j \le n/2$, $\ell \not\:\mid\; \!j$} \right\}
	\qquad\qquad\qquad\qquad\\\qquad\qquad\qquad\qquad
	\cup \left\{
	(1\:2\ldots 2k\ell-1)(2k\ell\;\,2k\ell+1) : \hbox{$1 < 2 k\ell < n$}
	 \right\}.\end{split}
\]
Note that $Z_\ell(n)$ consists of odd $\ell^\prime$-elements.
We can now apply Proposition~\ref{Prop:ThmAProp} to give
results about ordinary characters.

\begin{corollary}\label{Cor:Prop} Let $\lambda$ 
be a partition of $n$
and let $\ell > 1 $ be an odd natural number. Each of
the following conditions implies that $\lambda = \lambda'$:

\begin{indentedthmlist}
\item $\chi^\lambda$ vanishes on all cycles of even length
in $S_n$;

\item $\chi^\lambda$ vanishes on every element of $Z_\ell(n)$;

\item $\chi^\lambda$ vanishes on every odd $\ell^\prime$-element 
in $S_n$.
\end{indentedthmlist}
\end{corollary}

\begin{proof} (i) By hypothesis, $\omega_\lambda(s_{2i}) = 0$ for all $i$
such that $1 \le i \le n/2$. Hence, by
Proposition \ref{Prop:ThmAProp}(i), $\omega_\lambda(s_\mu) = 0$
whenever $\mu$ labels a conjugacy class of odd permutations. This implies that
$\chi^\lambda$ vanishes on every odd element in $S_n$, and so
$\chi^\lambda = \chi^\lambda \times \sgn$. Since
$\chi^\lambda \times \sgn = \chi^{\lambda^\prime}$
(see for instance
\cite[\S 6.6]{James}) we may deduce that~$\lambda = \lambda^\prime$.

(ii) This follows in the same way as (i), this time using
Proposition~\ref{Prop:ThmAProp}(iii).

(iii) This is merely a weaker version of the previous part.
\end{proof}

\begin{corollary}\label{Cor2:Prop}
Let $\lambda$ and $\mu$ be partitions of $n$.
Let $\ell > 2$ be a natural number. 
If $\chi^\lambda$ and $\chi^\mu$ agree on all $\ell'$-elements of $S_n$
then $\lambda = \mu$.
\end{corollary}

\begin{proof} By hypothesis
the central characters $\omega_\lambda$ and $\omega_\mu$ agree 
on the set $X_\ell(n)$ of generators of  $Z(\Q S_n)$ given in 
Proposition~\ref{Prop:ThmAProp}(ii). Hence $\omega_\lambda = \omega_\mu$ 
and so~$\lambda = \mu$.
\end{proof}

\subsection{Generating sets for $Z(\Q A_n)$}
Recall that the only conjugacy classes of $S_n$ which split in $A_n$
are those labelled by partitions of $n$ with odd distinct parts.
If $\mu$ is such a partition, let $s_\mu^+ \in Z(\Q A_n)$ and 
$s_\mu^- \in Z(\Q A_n)$ be the sums of the elements in the two 
associated conjugacy classes of $A_n$. 

For part (ii) of Theorem~\ref{Thm:DecDistinct} we need a generating
set for~$Z(\Q A_n)$ involving only conjugacy class sums labelled 
by $2^\prime$-permutations (that is,~permutations of odd order). 
Fortunately for us, the split classes consist of $2'$-elements, so 
they do not create any additional difficulties. Let
\[ \begin{split} X(n) = \bigl\{ s_\lambda : \lambda \vdash n,\, 
\text{all parts of $\lambda$ are odd} \bigr\} 
	\qquad\qquad\qquad\qquad\\\qquad\qquad\qquad\qquad\cup\,
\bigl\{ s_\lambda^+ : \lambda \vdash n,\, \text{$\lambda$ has odd distinct parts} \bigr\}.
\end{split}\]
We shall prove that $X(n)$ is a generating set for 
$Z(\Q A_n)$ as a $\Q$-algebra. To do this we need the following
lemma.

\begin{lemma}\label{Lemma:CycleProds}
Let $k,l \ge 3$ be natural numbers such that $k + l \le n+2$.
Define coefficients $c_\mu$ for $\mu \vdash n$ by
\[ s_k s_l = \sum_{\mu \vdash n} c_\mu s_\mu.\]
If $a,b \ge 2$ are natural numbers such that $a+b =k+l-2$ then
\[ c_{(a,b,1^{n-a-b})} = ab \min (k-1, l-1, a , b). \]
The only other conjugacy class sums $s_\mu$ with $\supp \mu \ge k+l-2$
which may appear as summands of $s_k s_l$ are $s_{(k+l-1)}$, which appears only if $k+l-1 \le n$,
and $s_{(k,l)}$, which appears only if $k+l \le n$.
\end{lemma}

The proof of this lemma is postponed to \S 2.3. 

\begin{proposition}\label{Prop:Prop2}
The elements of $X(n)$ generate $Z(\Q A_n)$ as a $\Q$-algebra.
\end{proposition}

\begin{proof}
It is sufficient to prove that if $\mu$ is a partition of $n$ with evenly 
many even parts and at least two even parts, then $s_\mu$ is in the 
$\Q$-algebra generated by~$X(n)$. 
Suppose that $\supp \mu = m$. By induction we may assume that all 
conjugacy class sums labelled by partitions with support at 
most~$m-1$ can be written as polynomials in elements of $X(n)$.

The hardest case occurs when $\mu$ has just two even parts
and all its other parts are of size $1$. Suppose that $m = 2r$.
Let $t = \left\lfloor r/2 \right\rfloor$ be
the number of partitions of~$2r$ into two even parts. 
By Lemma~\ref{Lemma:CycleProds}, 
if $1 \le i \le t$, then 
\[ \mfrac{1}{8} s_{2r-2i+1} s_{2i+1}  = \sum_{j=1}^i (r-j)j^2 s_{(2r-2j,2j)}
+ i\sum_{j=i+1}^t (r-j)j s_{(2r-2j,2j)} + y \]
where $y$ is a sum of conjugacy class sums labelled
by partitions which either have only odd parts, or 
have support at most $2r-1$. Therefore our inductive hypothesis
implies that for $1 \le i \le t$,
\[ u_i = \sum_{j=1}^i (r-j)j^2 s_{(2r-2j,2j)}
+ i \sum_{j=i+1}^t (r-j)j s_{(2r-2j,2j)} \]
is in the $\Q$-algebra generated by $X(n)$.
Now
\[ u_t - u_{t-1} = (r-t)t s_{(2r-2t,2t)} \]
and if $i < t$ then
\[ u_i - u_{i-1} = (r-i)i s_{(2r-2i,2i)} + \sum_{j=i+1}^t (r-j)j s_{(2r-2j,2j)}.
\]
Hence, by starting at $i=t$ and working down to $i=1$, we may
express each conjugacy class sum $s_{(2r-2i,2i)}$ as a polynomial 
in the elements of~$X(n)$. In particular, this shows that $s_\mu$ lies
in the $\Q$-algebra generated by~$X(n)$.

The other possibility is that~$\mu$ has two even parts,~$2u$ and~$2v$ say, 
and some further parts, not all of size $1$. Let~$\nu$ be the 
partition of $n-2u-2v$ obtained by removing the parts of size $2u$ 
and $2v$ from $\mu$. By induction $s_\nu$ and $s_{(2u,2v)}$ are in
the $\Q$-algebra generated by $X(n)$ and evidently
\[ s_\nu s_{(2u,2v)} = s_\mu + y \]
where~$y$ is a weighted sum of conjugacy class sums 
$s_\lambda$ for partitions $\lambda$ of support at most $m - 1$.
Hence~$s_\mu$ is in the $\Q$-algebra generated by $X(n)$.
This completes the inductive step.
\end{proof}

To obtain the expected corollary concerning ordinary characters of $S_n$
we need a small result about how characters of $S_n$ 
restrict to $A_n$. 

\begin{lemma}\label{Lemma:stupidlemma}
Let $\lambda$ and $\mu$ be  partitions of $n$. The restricted 
characters~$\chi^\lambda\res_{A_n}$ and $\chi^\mu\res_{A_n}$
agree if and only if $\lambda = \mu$ or $\lambda = \mu^\prime$.
\end{lemma}

\begin{proof} This may be proved using the ideas at the start of \S 4. Alternatively
see~\cite[\S 5.1]{FH} for a proof using only the language of character theory.
\end{proof}

\begin{corollary}\label{Cor:Prop2}
Let $\lambda$ and $\mu$ be partitions of $n$. If $\chi^\lambda$ and $\chi^\mu$ 
agree on all $2^\prime$-elements of $S_n$ then 
either $\lambda = \mu$ or $\lambda = \mu^\prime$.
\end{corollary}

\begin{proof}
By Proposition \ref{Prop:Prop2} and our usual argument with central
characters,
the hypothesis implies
that $\chi^\lambda(g) = \chi^\mu(g)$ for all~$g\in A_n$.
Now apply the previous lemma.
\end{proof}

\subsection{Proof of Lemma \ref{Lemma:CycleProds}}
We may assume that $k \ge l$. To find $s_ks_l$ we
shall first calculate the product $(12\ldots k)s_l$. As we are 
mainly interested in terms in this product whose support is 
exactly~$k + l - 2$, we start by looking at those~$l$-cycles 
which move exactly two members of $[1..k]$.
Let $i,j \in [1..k]$ be distinct numbers and let
\[ \tau = (i \: p_1 \ldots p_r \: j \: q_1 \ldots q_{l-2-r} )\]
where $0 \le r \le l-2$ and~$p_1,\ldots,p_r, q_1,\ldots, q_{l-2-r} \in 
[k+1.. n]$ are distinct numbers. Computation shows that if $i > 1$ then
\[ 
(1\:2\ldots k)\tau = 
(1 \ldots i-1 \: p_1 \ldots p_r \: j \ldots k) 
(i \: i+1 \ldots j-1 \: q_1 \ldots q_{l-2-r}),
\]
which is a product of cycles of lengths $(k+r) - (j-i)$ and $(j-i) + (l-2-r)$. 
If $i=1$ then we have 
\[ 
(1\:2\ldots k)\tau =
(1 \ldots j-1 \: q_1 \ldots q_{l-2-r}) 
   (j \: j+1 \ldots k \: p_1 \ldots p_r),
\]
so while there are some small differences,
the cycle structure of the product is unaltered.
It follows that if $2 \le a < (k+l-2)/2$ then
the number of ways to choose~$i$ and~$j$ so
that $(1\:2\ldots k)\tau$ has cycle type $(k+l-2-a,a)$ is
\begin{equation*}
\left\{ \begin{array}{@{\hspace{-0pt}}l@{\;\; :\;\;}l} a-r & a > r\\
								0	& a \le r \end{array} \right.
								+ \;
\left\{ \begin{array}{@{\hspace{-0pt}}l@{\;\; :\;\;}l} k-a+(l-2)-r & a > l-2-r \\
								0	& a \le l-2-r \end{array} \right.
\end{equation*}
where the first term comes from the case $(j-i)+(l-2-r) = k+l-2-a$ and the second from
the case $(j-i)+(l-2-r)=a$. 
This leaves out the possibility that
$2a  = k + l -2$. Then these two cases coincide, and the correct
number of choices for~$i$ and~$j$ is $a- r = (k+l-2)/2 - r$.

We now let $r$ vary between $0$ and $l-2$ and add up
the total number of choices for $i$ and $j$ so that
$(1\:2\ldots k)\tau$ has cycle type $(k+l-2-a,a)$.

\begin{indentedrmklist}
\item If $a \le l-2$ then the total number of choices for $i$ and $j$ 
is
\[ \sum_{0 \le r < a} (a-r) + 
\sum_{l-2 -a < r \le l-2} (k-a+(l-2)-r) = ak.\]

\item If $l-2 < a < (k+l-2)/2$ then the conditions needed
to get a positive contribution always hold in both cases, and the
total number of
choices is 
\[ \sum_{0 \le r \le m} \bigl( (a-r) + k-a + (l-2) -r \bigr) = (l-1)k. \] 

\item If $2a = k+l-2$ then 
the sum in the previous case double counts
every choice and so the total number of choices is $(l-1)k/2$.
\end{indentedrmklist}

We must also choose $p_1, \ldots, p_r$ and $q_1, \ldots, q_{m-2-r}$. 
Whatever the value of $r$, this can always
be done in exactly
$(n-k)^{\underline{l-2}}$ ways; here the notation~$x^{\underline{a}}$
stands for $x(x-1) \ldots (x-a+1)$ for  $x,a \in \N$.
Hence the total number of elements of cycle type 
$(k+l-2-a,a)$ in the product $(12\ldots k)s_l$ is
\[ (n-k)^{\underline{l-2}} \times
\left\{ \begin{array}{@{}l@{\;\; :\;\;}l} ak & 2 \le a \le l-2 \\
								 			(l-1)k & l-2 < a < (k+l-2)/2 \\
											(l-1)k/2 & a = (k+l-2)/2. 
											\end{array}\right. \]
To obtain the coefficient $c_{(k+l-2-a,a,1^{n-k-l+2})}$ 
we must multiply by $K_{(k,1^{n-k})}$ and 
then divide
by $K_{(k+l-2-a,a,1^{n-k-l+2})}$. This gives
\begin{align*} 
c_{(k+l-2,a,1^{n-k-l+2})} & = \min(a, l-1)\frac{ k(n-k)^{\underline{l-2}} 
	n^{\underline{k}}/k}{n^{\underline{k+l-2}}/(k+l-2-a)a} \\
	& = \min(a,l-1) (k+l-2-a)a
\end{align*}
as required.

It is easy to see that if $\tau$
is an $l$-cycle such that $(1\:2 \ldots k)\tau$ has support
strictly more than $k+l-2$ then either $\tau$ moves just one
element of $[1..k]$, in which case $(1\:2 \ldots k)\tau$
has cycle type $(k+l-1,1^{n-k-l+1})$, or $\tau$ moves no elements of $[1..k]$,
in which case $(1\:2 \ldots k)\tau$
has cycle type $(k,l,1^{n-k-l})$. This gives the final statement
in the lemma.\hfill$\Box$

\subsection{Some related problems}
We now pose some problems suggested by 
Corollaries~\ref{Cor:Prop},~\ref{Cor2:Prop} and~\ref{Cor:Prop2}. 
The reader keen to see the applications to the modular theory should skip
to \S 3.

\begin{problem}
What is the size $c_n$ of the smallest set $X \subseteq S_n$ 
such that if~$\chi^\lambda$ vanishes on all elements
of $X$ then $\lambda$ is self-conjugate?
\end{problem}

It follows from Corollary~\ref{Cor:Prop}(i) that~$c_n \le n/2$.
However, as Suzuki points out in \cite{Suzuki}, for~$n \le 14$, it suffices
to take $X = \left\{(12)\right\}$, so this result is not always
the best possible. An exhaustive search using the computer 
algebra package {\sc magma} shows that for $n \le 59$ one may 
take $X = \left\{(12), (1234)\right\}$, hence~$c_n = 1$ 
if $n \le 14$ and $c_n \le 2$ if $n \le 59$. (For $n=60$ 
there is a non-self-conjugate partition~$\lambda$ such that 
$\chi^\lambda((12)) = \chi^\lambda((1234)) = 0$, so it seems 
likely that $c_{60} = 3$.) We would ask for
an asymptotic formula for $c_n$, as its precise
behaviour may be too erratic to be easily described.

This problem may of course be posed with other properties in place of the condition that~$\lambda$ be self-conjugate. 
For example, one might ask  instead that~$\lambda$ 
be a $p$-core for a given prime $p$. Also one may restrict the possible 
set~$X$, for example by insisting that $X$ consist of $p^\prime$-elements
for a given prime~$p$. 

\begin{problem} 
What is the size $b_n$ of the smallest set 
$X \subseteq S_n$ such that if~$\chi^\lambda$ and~$\chi^\mu$ agree on all elements of~$X$ then $\lambda = \mu$? 
\end{problem}

Here Kramer's result  shows that $b_n \le n$, but
again this is not always the best possible.

\section{Consequences for decomposition matrices}
We are ready to prove Theorem~\ref{Thm:DecDistinct}. Let~$\phi_\nu$ 
be the Brauer character of the irreducible 
module~$D^\nu$ defined over a field
of prime characteristic~$p$ (for an introduction to Brauer characters 
see \cite[\S 2]{Navarro}). Suppose that in the decomposition 
matrix $D_p(n)$ of $S_n$ modulo~$p$, the rows labelled by 
partitions~$\lambda$ and~$\mu$ are equal. Adding up irreducible 
Brauer characters we find that if~$g$ is a~$p'$-element
of~$S_n$ then
\[ \chi^\lambda(g) = \sum_{\nu} D_p(n)_{\lambda\nu} \phi_\nu(g) 
= \sum_{\nu} D_p(n)_{\mu\nu} \phi_\nu(g) = \chi^\mu(g). \]
Thus $\chi^\lambda$ and $\chi^\mu$ agree on all $p'$-elements of $S_n$.
If~$p$ is odd then Corollary~\ref{Cor2:Prop} implies that $\lambda = \mu$. 
If $p=2$ then Corollary~\ref{Cor:Prop2} implies
that either~$\lambda = \mu$ or~$\lambda = \mu^\prime$. This completes
the proof.

We now generalise Theorem~\ref{Thm:DecDistinct} 
to Hecke algebras of symmetric groups.
Let~$F$ be any field (maybe of characteristic zero), and let
$q$ be an invertible element of $F$. Let $\mathcal{H}_{F,q}(S_n)$
be the corresponding Hecke algebra, as defined in~\cite[\S1.2]{MathasHecke} 
by deforming the group algebra $FS_n$. We may assume
that $q$ is a root of unity, as if not, 
$\mathcal{H}_{F,q}(S_n)$ is semisimple, and
so the analogue of Theorem~\ref{Thm:DecDistinct} is trivial.
Moreover, if $q=1$ then $\mathcal{H}_{F,q}(S_n) = FS_n$, so we may
also exclude this case. From now on we write $\mathcal{H}$ for
$\mathcal{H}_{F,q}(S_n)$.

Let $\ell$ be minimal such that $1+q+\ldots+q^{\ell-1} = 0$. 
The simple $\mathcal{H}$-modules are indexed by the $\ell$-regular
partitions of $n$ (and there are the expected analogues of
Specht modules), so the decomposition matrix of $\mathcal{H}$ 
has rows labelled by all partitions of $n$, and columns labelled
by the $\ell$-regular partitions of $n$. We need the following lemma,
which follows from the remarks just after
Proposition 2.6 in \cite{Richards}. 

\begin{lemma}
The $\Z$-span of the columns of the decomposition matrix of 
$\mathcal{H}$ is equal to the $\Z$-span of the 
columns of the ordinary character table of $S_n$ labelled by
the $\ell$-regular partitions of $n$. 
\hfill$\Box$
\end{lemma}

We can now prove the following analogue of Theorem~\ref{Thm:DecDistinct}.

\begin{theorem}\label{Thm:DecDistinctHecke}
Let $\mathcal{H}$ and $\ell$ be as above.

\begin{indentedthmlist}
\item If $\ell > 2$ then the rows of the decomposition matrix of  $\mathcal{H}$
are distinct.

\item If $\ell=2$ then the rows labelled by $\lambda$ and $\mu$ are the same
if and only if $\lambda = \mu$ or $\lambda = \mu^\prime$.
\end{indentedthmlist}
\end{theorem}

\begin{proof}
The previous lemma implies that
the rows of the decomposition matrix of $\mathcal{H}$ labelled
by partitions $\lambda$ and $\mu$ are equal if and only if
$\chi^\lambda(g) = \chi^\mu(g)$ for all $\ell'$-elements $g \in S_n$.
The result now follows from Corollary~\ref{Cor2:Prop}
and Corollary \ref{Cor:Prop2} in the same way as Theorem~\ref{Thm:DecDistinct}.
\end{proof}

We now turn to alternating groups. In odd characteristic the situation
appears to be quite difficult, and the obvious analogue of 
Theorem~\ref{Thm:DecDistinct} is false. There is however one result 
we can prove without any more work.

\begin{theorem}\label{Thm:Alt}
Let $n \in \N$. The rows of the decomposition
matrix of $A_n$ in characteristic $2$ are distinct.
\end{theorem}

\begin{proof}
Suppose that the rows labelled by the ordinary characters
$\chi$ and $\psi$ are equal. Then $\chi(g) = \psi(g)$ for all
$2'$-elements of $A_n$ and so by Proposition~2.5, $\chi(g) = \psi(g)$
for all $g \in A_n$. Hence $\chi = \psi$.
\end{proof}

\section{Specht modules and the alternating group}
In this section we prove Theorem~\ref{Thm:SpechtRes}.
Recall that this theorem states
that if $\lambda$ is a partition
of~$n$ and~$F$  
is an algebraically closed field of odd characteristic, then 
the Specht module~$S^\lambda$ 
is decomposable on restriction to~$A_n$ if and only if $S^\lambda$ 
is simple and $\lambda = \lambda^\prime$.

We begin with some Clifford theory. Let $F$ have characteristic $p$.
As~$p \not= 2$ the Specht module $S^\lambda$ is indecomposable
(see~\cite[Corollary 13.18]{James}).
Also, when $p \not=2$, the Sylow $p$-subgroups of $S_n$ are contained
in the alternating group $A_n$, so by Higman's
criterion (see \cite[Proposition 3.6.4]{Benson}),~$S^\lambda$ 
is relatively $A_n$-projective. 
Thus there exists an indecomposable $FA_n$-module~$U$ 
such that~$S^\lambda$ is a direct summand of the induced module
$U\ind_{A_n}^{S_n}$. We denote this by writing
\[ S^\lambda \:\bigl|\: U \Ind_{A_n}^{S_n}. \]
By Mackey's Lemma (see \cite[Theorem 3.3.4]{Benson}),
\begin{equation}
S^\lambda \Res_{A_n} \:\bigl|\: U \oplus U^t \label{Eq:Clifford}
\end{equation}
where $t$ is any odd element in $S_n$ and $U^t$ is the $A_n$-module
with the same underlying vector space as $U$, but with the action defined by 
$u\cdot g = ug^t$ for $u\in U^t$, $g\in A_n$.

We can now prove the `if' part of Theorem~\ref{Thm:SpechtRes}.
As $\lambda$ is self-conjugate, the restricted ordinary 
character~$\chi^\lambda\!\res_{A_n}$ splits as a sum of two ordinary
irreducible characters of $A_n$. Hence the 
Brauer character of $S^\lambda\!\res_{A_n}$ has at least two 
irreducible summands. (Notice that we have used that $F$ is sufficiently 
large here.) Furthermore, by Clifford's theorem on the restriction
of simple modules to normal subgroups, $S^\lambda \res_{A_n}$ is semisimple. 
Hence the restriction of $S^\lambda$ to $A_n$ has at least two non-trivial direct
summands. This is sufficient to prove the result.

It is however not hard to give a little more information. 
By~\eqref{Eq:Clifford}, $S^\lambda \res_{A_n}$ 
has at most two non-trivial direct summands. Therefore
\[ S^\lambda\Res_{A_n} = U \oplus U^t \]
and the Brauer characters of the simple summands $U$ and $U^t$ 
are the reduction modulo~$p$ of the two ordinary 
irreducible~$A_n$ characters associated to~$\lambda$.

We now turn to the `only if' part of Theorem~\ref{Thm:SpechtRes}.
Suppose that $S^\lambda \res_{A_n}$ is decomposable. Using
\eqref{Eq:Clifford}, we have 
\[ S^\lambda \otimes \sgn = U \Ind_{A_n}^{S_n} \otimes \sgn \cong 
\left( U \otimes \sgn\!\Res_{A_n} \right)\Ind_{A_n}^{S_n} \cong U \Ind_{A_n}^{S_n} \cong S^\lambda. \]
It follows that the ordinary character of $S^\lambda$ vanishes on
all odd $p'$-elements of $S_n$. By Corollary \ref{Cor:Prop}(iii)
this implies that $\lambda = \lambda^\prime$.
%(Alternatively one
%can note that the Brauer characters $S^\lambda$ and $S^{\lambda^\prime}$
%are equal and use Corollary \ref{Cor2:Prop}.)
Finally, since it is always the case that
\[ S^\lambda \otimes \sgn \cong (S^{\lambda^\prime})^\star, \]
where $(S^{\lambda^\prime})^\star$ is the dual module to $S^{\lambda^\prime}$ 
(see \cite[Theorem 8.15]{James}), $S^\lambda$
is self-dual, and hence simple by Theorem~23.1 of \cite{James}.
This completes the proof of Theorem~\ref{Thm:SpechtRes}.

\subsection{Theorem \ref{Thm:SpechtRes} for non-algebraically closed fields}
\label{sect:nonalgclosed}
We now consider Specht modules
defined over an arbitrary field $F$ of odd characteristic $p$. 
As explained in Remark 2 of \S 1.2, we shall use a result due to Fayers
on the irreducibility of Specht modules. To state it we 
need two final pieces of notation. If $\lambda$
is a partition, and $\alpha$ is a node of  the diagram of $\lambda$, let
$h_{\alpha}$ be the hook-length of the hook on $\alpha$.
(See \cite{FayersIf} or \cite[Chapter~18]{James} for the definition of 
hooks in partitions.) Given $n\in \N$, let~$(n)_p = p^a$ if $p^a$
is the highest power of~$p$ which divides~$n\in \N$. 

\begin{theorem}[Fayers]\label{Thm:Fayers}%NOTE needs (Fayers) in LMS class
Let $F$ be a field of odd characteristic. The Specht module~$S^\lambda$ 
defined over $F$ is reducible if and only if~$\lambda$ contains a 
node~$\alpha$, a node~$\beta$ in the same row as~$\alpha$,
and a node~$\gamma$ in the same column as~$\alpha$, such 
that $p \left|\right. h_\alpha$, $(h_\alpha)_p \not= (h_\beta)_p$ and
$(h_\alpha)_p \not= (h_\gamma)_p$.\hfill$\Box$
\end{theorem}

Using the `if' direction of this theorem (proved in \cite{FayersIf})
we may prove the following lemma.

\begin{lemma}\label{Lemma:JMlemma}
Let $\lambda = (\lambda_1,\ldots,\lambda_k)$ be a self-conjugate 
partition of $n$. Suppose that there is a node on the main 
diagonal of~$\lambda$ whose hook-length is divisible by~$p$. 
Then $S^\lambda$ is reducible.
\end{lemma}

\begin{proof} For $1\le i \le k$ and $1 \le j \le \lambda_i$ 
let $h_{ij}$ be the hook-length
of the node in position~$(i,j)$ of $\lambda$. Suppose that
$(h_{ss})_p = p^c$ where $c \ge 1$. 
By Theorem~\ref{Thm:Fayers},
if $S^\lambda$ is irreducible
then every node in row $s$ and every node in column $s$ has
hook length exactly divisible by $p^c$. By considering the 
subdiagram of~$\lambda$
obtained by taking all nodes in positions $(a,b)$ for $a, b \ge s$ 
we obtain a partition $\mu$ whose hook lengths satisfy
\[ (h_\alpha)_p = p^c \]
for all nodes~$\alpha$ in the first row and column of $\mu$. Let~$\mu_1 = l$.
As $p^c \mid h_{1l}$ we must have $\mu_l^\prime > 1$. Hence $\mu_1 = \mu_2$,
and so $h_{21} = h_{11} - 1$. But both these
hook lengths are supposed to be divisible by~$p^c$, so we have
reached a contradiction.
\end{proof}

We are now ready to prove the following generalisation 
of Theorem~\ref{Thm:SpechtRes}.

\begin{theorem}\label{Thm:SmallField} Let $F$ be a field
of odd characteristic. Let $\lambda$ be a partition of~$n \in \N$
with main diagonal hook lengths $q_1, \ldots, q_r$.
The~$FS_n$-Specht module~$S^\lambda$ is decomposable
on restriction to~$A_n$ if and only if all of the following conditions hold:

\begin{indentedthmlist}
\item $\lambda$ is self-conjugate,
\item $(-1)^{(n-r)/2}q_1 \ldots q_r$ has a square root in $F$,
\item $S^\lambda$ is simple (we will  see this implies
that $p$ does not divide any of the~$q_i$).
\end{indentedthmlist}
\end{theorem}

\begin{proof} 
By the proof of Theorem \ref{Thm:SpechtRes}, if $S^\lambda\res_{A_n}$
is decomposable then~(i) and~(iii) hold.
Let $\chi_1$ and $\chi_2$ be the two ordinary
characters of $A_n$ associated to $\lambda$: at $p^\prime$ elements
these are the Brauer characters of the summands of~$S^\lambda\res_{A_n}$.
Proposition 5.3 of \cite{FH}
tells us that
\[ \mfrac{1}{2} (-1)^{(n-r)/2} \pm \mfrac{1}{2} 
\sqrt{(-1)^{(n-r)/2}q_1 \ldots q_r} \]
are the values of $\chi_1$ and $\chi_2$ on elements of the 
conjugacy class labelled by~$(q_1, \ldots, q_r)$. 
If any of the $q_i$ are divisible by $p$ then Lemma \ref{Lemma:JMlemma}
implies that~$S^\lambda$ is reducible, a contradiction. Hence 
this conjugacy class is $p$-regular. Therefore
\[ (-1)^{(n-r)/2}q_1 \ldots q_r \] 
has a square root in $F$, which gives (ii).

Conversely if all the conditions hold then the proof of the `if' part
of Theorem \ref{Thm:SpechtRes} shows that $S^\lambda \res_{A_n}$ decomposes.
Where before we used that $F$ was sufficiently large, now we merely
use the fact that if $\chi$ is the character, in the na{\"i}ve sense, of an
irreducible representation of a group $G$ over a field~$E$ of prime
characteristic, then the representation can be defined over a subfield~$F$ 
of~$E$ if and only if the values of~$\chi$ lie in~$F$. 
(For a proof of this statement see \cite[Theorem~9.14]{IsaacsChars}.)
\end{proof}

We conclude by noting that the last
result has an especially easy form for Specht modules
labelled by hook partitions.

\begin{corollary}
Let $F$ be a field of odd prime characteristic, and let
$1 < r < n-1$. The Specht module $S^{(n-r,1^r)}$ decomposes
on restriction to~$A_n$ if and only if $n = 2r+1$,~$p$ does
not divide~$n$, and~$(-1)^{(n-1)/2}\thinspace n$ has a square root in~$F$.\hfill$\Box$
\end{corollary}

\section{Acknowledgements}
%\begin{acknowledgements} %Used for LMS class
I should like to thank
Matt Fayers for his comments
on Theorem \ref{Thm:SpechtRes} and for alerting me at an early stage
to his proof of the James--Mathas conjecture.
%\end{acknowledgements}

\def\cprime{$'$} \def\Dbar{\leavevmode\lower.6ex\hbox to 0pt{\hskip-.23ex
  \accent"16\hss}D}

\end{document}